\documentclass[
  12pt,
  paper=a4,
  numbers=noenddot,
  headings=normal
]{scrartcl}

\usepackage[T1]{fontenc}
\usepackage{libertinus}
\usepackage{microtype}

\usepackage{amsmath,amssymb,amsthm,mathtools}

\usepackage[utf8]{inputenc}
\usepackage{enumerate}
\usepackage{booktabs,tabularx}
\usepackage{xcolor}
\usepackage{stmaryrd}
\usepackage[numbers,sort&compress]{natbib}
\usepackage{tikz-cd}
\usepackage{tikz}
\usepackage{float}
\usepackage{extarrows}

\usepackage{authblk}

\newcommand{\email}[1]{\href{mailto:#1}{#1}}

\newcommand{\longtwoheadrightarrow}[1]{%
  \mathrel{\overset{#1}{\relbar\joinrel\twoheadrightarrow}}%
}

\setkomafont{section}{\normalfont\Large\bfseries}
\setkomafont{subsection}{\normalfont\large\bfseries}
\setkomafont{subsubsection}{\normalfont\normalsize\bfseries}

\RedeclareSectionCommand[
  beforeskip=-3.5ex plus -1ex minus -.2ex,
  afterskip=1.7ex plus .2ex
]{section}

\RedeclareSectionCommand[
  beforeskip=-2.7ex plus -1ex minus -.2ex,
  afterskip=1.0ex plus .2ex
]{subsection}

\usepackage[hidelinks]{hyperref}
\usepackage[nameinlink,capitalize,noabbrev]{cleveref}
\hypersetup{
  pdftitle={Graphs with isomorphic maximal abelian covers},
  pdfauthor={Frederico Cancado, Gabriel Coutinho, Thomas Jung Spier}
}

\newtheorem{theorem}{Theorem}[section]
\newtheorem{lemma}[theorem]{Lemma}
\newtheorem{proposition}[theorem]{Proposition}
\newtheorem{corollary}[theorem]{Corollary}

\theoremstyle{definition}

\newtheorem{example}[theorem]{Example}

\newcommand{\Z}{\mathbb Z}
\newcommand{\R}{\mathbb R}

\newcommand{\ab}{\mathrm{ab}}
\newcommand{\Aut}{\operatorname{Aut}}
\newcommand{\Deck}{\operatorname{Deck}}
\newcommand{\Cov}{\operatorname{Cov}}
\newcommand{\BD}{\operatorname{BD}}
\newcommand{\im}{\operatorname{im}}
\newcommand{\rank}{\operatorname{rank}}
\newcommand{\Core}{\operatorname{Core}}

\newcommand{\one}{\mathbf 1}
\newcommand{\hml}[1]{\llbracket #1 \rrbracket}
\newcommand{\word}{\mathfrak w}

\newcommand{\ba}{\gamma}

\title{\LARGE Complete characterization of graphs with common maximal abelian cover}

\author[1]{Frederico Cançado\textsuperscript{*}}
\author[1]{Gabriel Coutinho}
\author[2]{Thomás Jung Spier}

\affil[1]{Department of Computer Science,
Federal University of Minas Gerais, Brazil}

\affil[2]{Department of Mathematics,
Federal University of Santa Catarina, Brazil}

\affil[ ]{\small
\email{fredericocancado@dcc.ufmg.br}\quad
\email{gabriel@dcc.ufmg.br}\quad
\email{thomas.j.s@ufsc.br}}

\affil[ ]{\small
\textsuperscript{*}Corresponding author}

\date{\today}

\begin{document}

\maketitle

\begin{abstract}
The maximal abelian cover of a connected graph is the quotient
of its universal cover by the commutator subgroup of its
fundamental group. We provide a complete classification of finite connected graphs with
isomorphic maximal abelian covers. For trees or for connected graphs with at least two cycles,
the cover determines the base up to isomorphism; in cycle
rank one, it determines the primitive cyclic pattern of rooted trees attached to the unique cycle. We also characterize when every
automorphism descends to the base and give an orbit-size
criterion for splitting the canonical lift extension, answering two questions raised by Baez in 2026.
\end{abstract}

\section{Introduction}\label{sec:classification}

How much of a finite graph is determined by its maximal abelian cover? The corresponding question for universal covers is well understood: nonisomorphic graphs can have the same universal covering tree, and Leighton's theorem identifies this relation with the existence of a common finite cover \cite[Section~3]{Leighton1982}. We show that the maximal abelian cover retains substantially more information. Our main result is that, apart from unicyclic graphs, its abstract isomorphism type determines the base graph.

Abelian covers have long been studied through homology and periodic graph realizations. Biggs \cite{Biggs1984} developed homological graph coverings, while Kotani and Sunada \cite{KotaniSunada2000,KotaniSunada2001} related abelian covers to graph tori and harmonic realizations; see also Sunada's monograph \cite{Sunada2013}. Baez \cite[Theorems~3 and~8]{Baez2026} gives an explicit cycle-space embedding of the maximal abelian cover of a bridgeless graph. Recent work of Li, Magee, Sabri, and Thomas \cite{LiEtAl2025} characterizes the eigenvalues of maximal abelian covers in terms of the finite base graph. Our question concerns reconstruction from the cover itself, rather than from a chosen realization or its spectrum.

These questions also arise in crystal chemistry, where
periodic nets model the connectivity of crystal frameworks.
Beukemann and Klee's work on maximal abelian covers
\cite{BeukemannKlee1992}, Eon's studies of Euclidean embeddings
\cite{Eon1999,Eon2011Embeddings}, and his joint work with
Moreira de Oliveira on automorphism groups
\cite{MoreiraDeOliveiraEon2011} provide important background
on the relationships between finite quotient graphs,
infinite periodic nets, and their geometric realizations.
This crystallographic perspective motivates the
reconstruction and symmetry questions considered here.

The main difficulty is to recover the covering fibres without
being given the projection. Lifting an automorphism of a base
graph is a classical problem \cite{MalnicEtAl2000}; for the
maximal abelian cover, every base automorphism lifts
\cite[Lemma~9]{Baez2026}. The converse is not automatic.
Baez asks when every automorphism of the maximal abelian cover
descends to the base \cite[Section~8]{Baez2026}.
For the reconstruction problem, an earlier statement for finite
connected bases of minimum degree at least three appears in
Beukemann and Klee's work on maximal abelian covers
\cite[Lemmas~1--2 and the first theorem on p.~42]{BeukemannKlee1992}.
For the descent question, an affirmative answer for finite
\(c\)-connected bases---that is, connected graphs in which every
vertex lies on a cycle or loop---of cycle rank at least two
follows from Eon's earlier work
\cite[proof of Theorem~6.1, pp.~62--63]{Eon2007}.
Eon's argument uses the fact that bounded-displacement
automorphisms send geodesic fibres to parallel fibres
\cite[Theorem~5.1]{Eon2007}.

The earlier reconstruction and descent arguments draw on
Euclidean embeddings and their translation structure
\cite{BeukemannKlee1992}, or on the geometry of geodesics
and fibres \cite{Eon2007}. Our proof combines integral
harmonic coordinates, an elementary maximum principle,
and the rank formula for finite graph coverings.
For finite connected graphs of cycle rank at least two,
this gives a self-contained proof of Theorem~\ref{thm:deck}
under the weaker hypothesis of minimum degree at least two.

A substantial part of our contribution lies in bringing together and
interpreting results scattered across graph theory and
crystallography, expressed in different languages and under
different hypotheses. We clarify their historical development,
contributions, and limitations, and combine these advances
with our results to give a unified account culminating in
a complete characterization of finite connected graphs
with isomorphic maximal abelian covers and full answers
to two questions posed by Baez
\cite[Section~8]{Baez2026}.

We use integral harmonic coordinates, an elementary maximum
principle, and the rank formula for finite graph coverings.
For finite connected graphs of cycle rank at least two,
Theorem~\ref{thm:deck} assumes only minimum degree at least two.
This condition is weaker than Eon's \(c\)-connectedness and the
minimum-degree-three condition on Beukemann--Klee's
finite quotient graphs.
Combining this result with an analysis of the rooted trees
attached to the core, we prove our main classification
theorem and give a complete answer to Baez's descent question
for finite connected graphs: every automorphism descends
precisely when the base is a tree or all its attached rooted
trees have trivial root-preserving automorphism groups.
In Section~\ref{sec:baez}, we also answer Baez's second question,
concerning the semidirect-product structure of the
covering-symmetry group.

\subsection{Notation and covering background}\label{subsec:background}

All base graphs are finite, undirected, connected, and nonempty. Covering
graphs may be infinite, but are connected and locally finite. Loops and
multiple edges are allowed. Each edge is represented by two distinct
reverse arcs $a$ and $\bar a$, with initial and terminal vertices
$\iota(a)$ and $\tau(a)$. We write
\[
  A_G^+(v)=\{a\in A(G):\iota(a)=v\}.
\]
Thus $\deg_G(v)=|A_G^+(v)|$, and a loop contributes twice to the degree.
Isomorphisms preserve vertices, arcs, incidence, and reversal. These
conventions agree with the directed-edge treatment in
\cite[Section~1.1.5]{GrossTucker1987}.

A \emph{graph map} $f:X\to G$ consists of maps
$f_V:V_X\to V_G$ and $f_A:A_X\to A_G$ satisfying
\[
  f_V\circ\iota_X=\iota_G\circ f_A,\qquad
  f_V\circ\tau_X=\tau_G\circ f_A,\qquad
  f_A(\bar a)=\overline{f_A(a)}
  \quad(a\in A_X).
\]
We write $f$ for both maps when their domains are clear.
A \emph{covering} $p:X\to G$ is a surjective graph map inducing a
bijection $A_X^+(x)\to A_G^+(p(x))$ at every vertex $x$. We say the vertex $x$ is above $v$, or is a lift of $v$ when $p(x)=v$. 

For a walk $P=a_1\cdots a_m$ in $G$ and a vertex $x$
above its initial vertex, a \emph{lift} of $P$ from $x$
is a walk $\widetilde P=\widetilde a_1\cdots\widetilde a_m$
in $X$ starting at $x$ and satisfying
$p(\widetilde a_j)=a_j$ for every $j$.
Such a lift exists and is unique. Different occurrences
of the same vertex in $P$ may lift to different vertices
of $X$; in particular, a closed walk need not lift to a
closed walk.

For the voltage-graph formulation of unique walk lifting, see
\cite[Theorems~2.1.1 and~2.4.1]{GrossTucker1987}. Lifting a path between
two base vertices gives a bijection between their fibres; in particular,
all vertex fibres of a connected covering have the same cardinality.

For a subgraph, we distinguish a lifted copy from the
lift of a walk traversing it. Let $H\subseteq G$ be a
nonempty connected subgraph, with inclusion
$i:H\hookrightarrow G$. A \emph{lift} of this inclusion
is a graph map $\widetilde i:H\to X$ such that
\[
  p\circ\widetilde i=i.
\]
We also call its image $\widetilde H=\widetilde i(H)$
a \emph{lift} of $H$. Equivalently,
$p|_{\widetilde H}:\widetilde H\to H$ is a graph
isomorphism. In particular, $\widetilde H$ is connected
and has the same numbers of vertices and edges as $H$.
The full preimage $p^{-1}(H)$ need not be connected or
have the same numbers of vertices and edges as $H$.

For $v\in V(H)$ and $x\in p^{-1}(v)$, a lift with
$\widetilde i(v)=x$ exists if and only if every closed walk
in $H$ based at $v$ lifts to a closed walk starting at $x$.
When it exists, it is unique
\cite[Propositions~1.33--1.34]{Hatcher2002}.
In particular, a tree has a unique lifted copy through
any prescribed vertex above one of its vertices.

The fundamental group $\pi_1(G,r)$ consists of closed walks at a fixed
root $r$, modulo insertion and deletion of immediate backtracks. Its
product is concatenation followed by reduction. A spanning tree gives
free generators indexed by the edges outside the tree, so its rank is
\[
  \beta(G)=|E_G|-|V_G|+1,
\]
the \emph{cycle rank} of $G$; see
\cite[Section~2.5.3]{GrossTucker1987} and
\cite[Proposition~1A.2]{Hatcher2002}. The universal cover $T_G\to G$
is the tree whose vertices are reduced walks starting at $r$, with
projection to their terminal vertices. Its deck transformation group
is naturally isomorphic to $\pi_1(G,r)$.

A \emph{deck transformation} of $p:X\to G$ is an automorphism $d$ of
$X$ with $p\circ d=p$; their group is denoted by $\Deck(p)$. A connected
cover is \emph{regular} when this group is transitive on each fibre.
More generally, if a group $\Gamma\leq\Aut(X)$ acts freely on vertices
and without edge inversions, the quotient $X/\Gamma$ has the vertex
and arc orbits as its vertices and arcs. The natural projection is a
regular covering. Conversely, a regular covering is recovered from
its deck-group orbits; see
\cite[Theorems~2.2.1--2.2.2]{GrossTucker1987}.

The covering--subgroup correspondence identifies connected covers of
$G$ with subgroups of $\pi_1(G,r)$, up to conjugacy. A normal subgroup
$N$ gives a regular cover with deck group $\pi_1(G,r)/N$
\cite[Theorem~1.38 and Proposition~1.39]{Hatcher2002}.
In particular, the \emph{maximal abelian cover} is
\begin{equation}\label{eq:cover}
  G^{\ab}=T_G/[\pi_1(G,r),\pi_1(G,r)],
\end{equation}
where $[\Gamma,\Gamma]$ denotes the subgroup generated by the
commutators $ghg^{-1}h^{-1}$. Its deck group is
\[
  \pi_1(G,r)^{\ab}
  :=\pi_1(G,r)/[\pi_1(G,r),\pi_1(G,r)]
  \cong H_1(G;\Z)\cong\Z^{\beta(G)},
\]
where $H_1(G,\Z)$ denotes the \emph{first homology group} of $G$ (see Section \ref{sec:walks}), and the first isomorphism is known as Hurewicz theorem \cite[Theorem~2A.1, p.~166]{Hatcher2002}. The word \emph{maximal} refers to the covering order: every connected
regular abelian cover of $G$ is covered by $G^{\ab}$ over $G$.
In Section~\ref{sec:walks} we give the equivalent walk model, in which
walks with the same endpoint are identified precisely when they have
the same signed edge counts; see \cite[Section~2 and Lemma~4]{Baez2026}.

For later comparison with voltage graphs, choose a spanning tree and
orient the $b:=\beta(G)$ remaining edges. Define the \emph{voltage function} $\alpha: A(G)\to \Z^b$ by assigning voltage $0\in\Z^b$
to each tree arc and the $j$th standard basis vector to the chosen orientation of the $j$th
cotree edge, with opposite voltage on its reverse. The derived
graph has vertices $(v,z)\in V_G\times\Z^b$ and arcs
\[
  (\iota(a),z)\longrightarrow(\tau(a),z+\alpha(a)).
\]
A closed walk has zero voltage exactly when its image in
$\pi_1(G,r)^{\ab}$ is zero, so this derived graph is $G^{\ab}$.
This is the spanning-tree normalization of ordinary voltages
\cite[Sections~2.1 and~2.5.2--2.5.3]{GrossTucker1987}.

\subsection{The classification and proof outline}\label{subsec:main-result}

For a unicyclic graph $G$, list the vertices of its unique cycle in
cyclic order as $v_0,\ldots,v_{m-1}$. A loop is a cycle of length one,
and two parallel edges form a cycle of length two. After deleting
the cycle edges, let $T_i$ be the remaining tree rooted at $v_i$.
Regard
\[
  w(G)=(T_0,\ldots,T_{m-1})
\]
as a word in rooted-tree isomorphism types. Its \emph{primitive
decoration word} is its shortest repeating block, considered up to
cyclic rotation and reversal. Denote this equivalence class by
$\word(G)$.

\begin{theorem}\label{thm:classification}
Let $G,H$ be finite connected graphs with the conventions above.
\begin{enumerate}[(a)]
\item If $\beta(G)=\beta(H)=0$, then $G^{\ab}\cong H^{\ab}$ if and only if
      $G\cong H$.
\item If $\beta(G)=\beta(H)=1$, then $G^{\ab}\cong H^{\ab}$ if and only if
      $\word(G)=\word(H)$.
\item If $\beta(G),\beta(H)\ge2$, then $G^{\ab}\cong H^{\ab}$ if and only if
      $G\cong H$.  \label{eq:rigidity}
\end{enumerate}
Graphs with different cycle ranks cannot have isomorphic maximal abelian
covers. In particular, nonisomorphic graphs can have isomorphic maximal
abelian covers only in the unicyclic case.
\end{theorem}

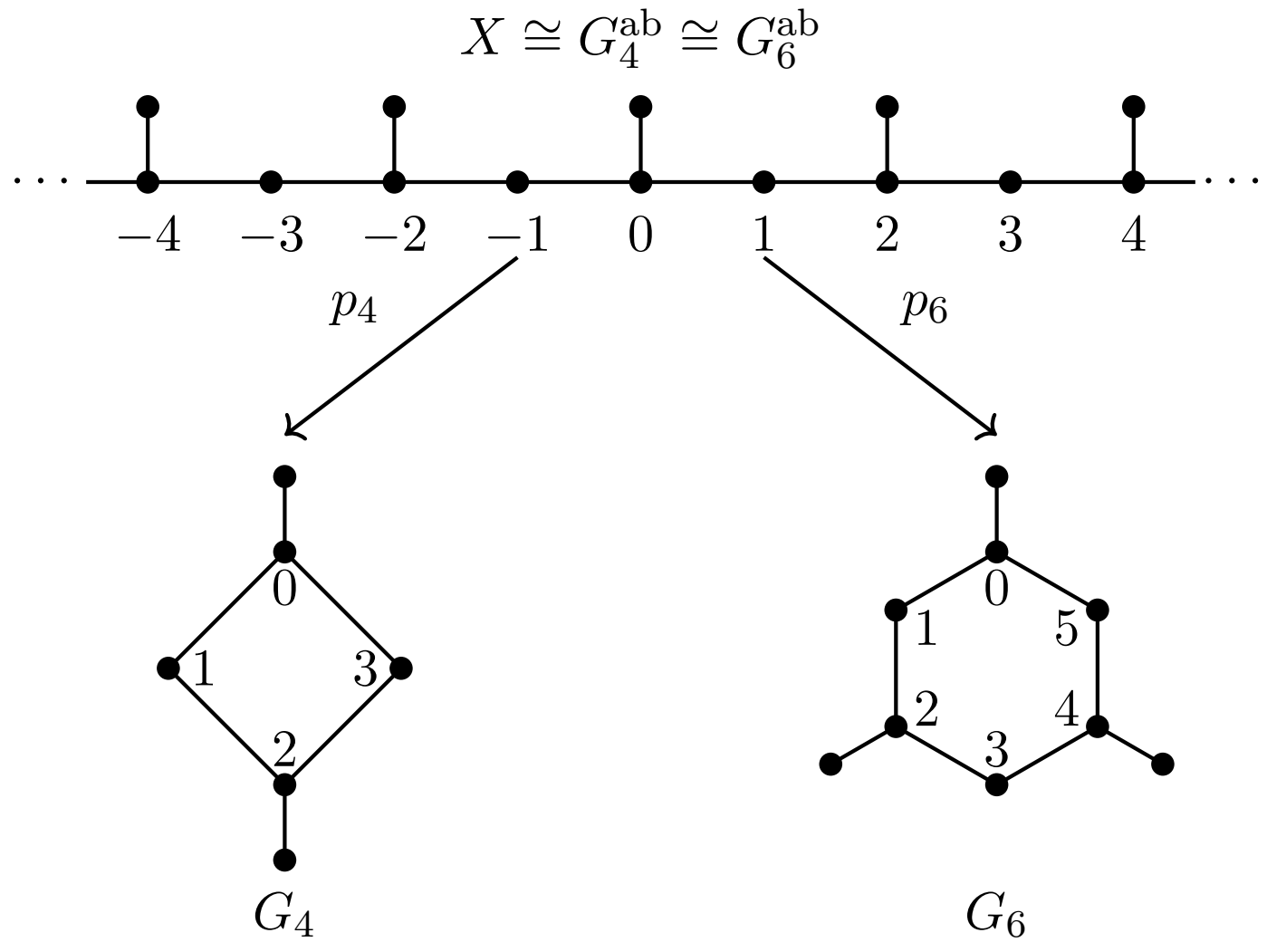
\begin{figure}[H]
\centering
\begin{tikzpicture}[
  vertex/.style={circle, fill=black, inner sep=1.7pt},
  thick, font=\small
]
  \node at (0,3.45) {$X\cong G_4^{\ab}\cong G_6^{\ab}$};
  \draw (-4.05,2.4)--(4.05,2.4);
  \foreach \n in {-4,...,4}{
    \node[vertex, label=below:{$\n$}] at ({0.9*\n},2.4) {};
  }
  \foreach \n in {-4,-2,0,2,4}{
    \draw ({0.9*\n},2.4)--++(0,0.55) node[vertex] {};
  }
  \node at (-4.35,2.4) {$\cdots$};
  \node at (4.35,2.4) {$\cdots$};

  \draw[->] (-0.9,1.85)--(-2.6,0.55)
    node[midway, above left] {$p_4$};
  \draw[->] (0.9,1.85)--(2.6,0.55)
    node[midway, above right] {$p_6$};

  \begin{scope}[shift={(-2.6,-1.15)}]
    \foreach \i/\angle in {0/90,1/180,2/270,3/0}{
      \coordinate (g\i) at (\angle:0.85cm);
      \node[vertex] at (g\i) {};
      \node at (\angle:0.59cm) {$\i$};
    }
    \draw (g0)--(g1)--(g2)--(g3)--cycle;
    \draw (g0)--++(90:0.55cm) node[vertex] {};
    \draw (g2)--++(270:0.55cm) node[vertex] {};
    \node at (0,-1.8) {$G_4$};
  \end{scope}

  \begin{scope}[shift={(2.6,-1.15)}]
    \foreach \i/\angle in {0/90,1/150,2/210,3/270,4/330,5/30}{
      \coordinate (h\i) at (\angle:0.85cm);
      \node[vertex] at (h\i) {};
      \node at (\angle:0.59cm) {$\i$};
    }
    \draw (h0)--(h1)--(h2)--(h3)--(h4)--(h5)--cycle;
    \foreach \i/\angle in {0/90,2/210,4/330}{
      \draw (h\i)--++(\angle:0.55cm) node[vertex] {};
    }
    \node at (0,-1.8) {$G_6$};
  \end{scope}
\end{tikzpicture}
\caption{Two nonisomorphic unicyclic graphs and their common maximal
abelian cover. On the path vertices, $p_m(n)=n\bmod m$; each leaf
maps to the leaf at the corresponding cycle vertex ($m=4,6$).}
\label{fig:unicyclic-same-abelian-cover}
\end{figure}

For example, let $G_4$ and $G_6$ be obtained by attaching one leaf at
every other vertex of a $4$-cycle and a $6$-cycle, respectively.
These graphs are not isomorphic, but both maximal abelian covers
are the two-way infinite path with one leaf attached at every
other vertex, as shown in Figure~\ref{fig:unicyclic-same-abelian-cover}.
Their primitive decoration words are the same.

In contrast, the triangular prism and $K_{3,3}$ have the same
universal cover but, by the theorem, different maximal abelian covers.

The central step is an intrinsic description of the deck group.
For a connected graph $X$, let $\BD(X)$ be the group of automorphisms
$\ba$ for which $d_X(x,\ba x)$ is bounded independently of $x$. Such
\emph{bounded automorphisms} have been studied by Trofimov
\cite{Trofimov1984} and Godsil et al.\ \cite{GodsilEtAl1989}.
We prove that, when $G$ has minimum degree at least two and
$\beta(G)\ge2$,
\[
  \BD(G^{\ab})=\Deck(G^{\ab}\to G).
\]
Since bounded displacement is preserved by every graph isomorphism,
this equality recovers the covering fibres and hence the base.

Our proof uses integral coordinates for the classical harmonic realization of Kotani and Sunada \cite[Section~4, equations~(28)--(29) and~(35)]{KotaniSunada2001}. Their fibres are the lifted bridge trees, reflecting the bridge-contraction property of the graph Abel--Jacobi map \cite[Lemma~4.6 and Corollary~4.7]{BakerNorine2007}. The additional step is to show that every bounded-displacement automorphism acts by a constant coordinate translation, and then to identify the resulting translation group with the canonical deck group. Integrality reduces the harmonic-function argument to an elementary maximum principle.

Section~\ref{sec:walks} develops the walk model and cycle coordinates.
Section~\ref{sec:coordinates} proves the properties of the harmonic
realization needed in the argument. Section~\ref{sec:deck} reconstructs
the deck group and proves descent of isomorphisms for graphs without
leaves. Section~\ref{sec:classification-proof} recovers the attached
trees and handles the small cycle ranks, completing the proof of
Theorem~\ref{thm:classification}.
Section~\ref{sec:baez} proves the descent and splitting criteria.
Section~\ref{sec:scope} discusses further questions and
distinguishes isomorphic maximal abelian covers from the
existence of a common connected regular abelian cover.

\section{Walks, cycles, and the canonical cover}\label{sec:walks}
The standard covering-space and spanning-tree foundations can be found in
\cite[Section~2.5.3]{GrossTucker1987} and
\cite[Section~1.3 and Proposition~1A.2]{Hatcher2002}. We give the
particular constructions needed here explicitly.

A spur in a walk is a consecutive pair $a\bar a$; a walk is reduced if it has
no spur. Combinatorially, $\pi_1(G,r)$ is the group of closed walks at $r$ modulo
insertion and deletion of spurs, with concatenation as product. For a group
$\Gamma$, its commutator subgroup $[\Gamma,\Gamma]$ is generated by the elements
$ghg^{-1}h^{-1}$, and its abelianization is $\Gamma^{\ab}=\Gamma/[\Gamma,\Gamma]$.

Choose one orientation of every edge of $G$. Let $C_1(G;\Z):=\Z^{E_G}$ and
$C_0(G;\Z):=\Z^{V_G}$, with boundary map
\[
  \partial:C_1(G;\Z)\longrightarrow C_0(G;\Z),\qquad
  (u\to v)\longmapsto\one_v-\one_u.
\]
Let $\epsilon:A(G)\to C_1(G;\Z)$ assign to each arc its signed unit
edge vector. Thus $\epsilon(\bar a)=-\epsilon(a)$, and
$\epsilon(a)=\one_e$ when $a$ agrees with the chosen orientation
of its underlying edge $e$. In the canonical bases, $\partial$ is
the incidence matrix, with $-1$ at the initial vertex and $+1$ at
the terminal vertex of an oriented edge.

For a walk $P=a_1\cdots a_s$, its \emph{signed edge-count vector},
or \emph{chain}, is
\[
  c(P):=\sum_{i=1}^{s}\epsilon(a_i).
\]
Thus a walk from $u$ to $v$ satisfies $\partial c(P)=\one_v-\one_u$. The integral homology group is
\[
  H_1(G;\Z):=\ker\partial.
\]
We use the standard Euclidean inner product on real edge vectors. For the
integral cycle and cut lattices, see \cite[Section 1]{BacherEtAl1997}; for the associated
graph tori and abelian covers, see \cite{KotaniSunada2000}.

The spanning-tree basis of the integral flow lattice, which is
$\ker\partial$ in our notation, is given in
\cite[Lemma~2, p.~182]{BacherEtAl1997}. The associated fundamental
closed walks are used in \cite[Section~2.5.3, p.~92]{GrossTucker1987}.
The following standard lemma records their compatibility with signed
edge counts.

\begin{lemma}[{\cite[Lemma~2]{BacherEtAl1997};
\cite[Proposition~1A.2 and Theorem~2A.1]{Hatcher2002}}]
\label{lem:tree-basis}
Fix a root $r$ and a spanning tree $S$. Let $P_v$ be the unique tree path
from $r$ to $v$. Orient the edges outside $S$ as $e_1,\ldots,e_b$, where
$b=\beta(G)$, and put
\[
  z_j=c(P_{\iota(e_j)})+\epsilon(e_j)-c(P_{\tau(e_j)}).
\]
Then the $z_j$ form a $\Z$-basis of $\ker\partial$, and every entry of each
$z_j$ belongs to $\{0,1,-1\}$. The fundamental group of $G$ is free on the
closed walks
$P_{\iota(e_j)}e_jP_{\tau(e_j)}^{-1}$, and signed edge counts induce an
isomorphism
\begin{equation}\label{eq:homology}
  \pi_1(G,r)^{\ab}\cong\ker\partial\cong\Z^b.
\end{equation}
\end{lemma}

Baez's \emph{walk model} identifies homology classes with signed edge-count
vectors. Applying the covering--subgroup correspondence to the commutator
subgroup gives regularity and the deck group below. The simplicity
assertion follows directly from the fact that distinct outgoing arcs
change the chain coordinate by distinct nonzero signed unit vectors. We present the walk model by the following proposition:

\begin{proposition}[{\cite[Section~2 and Lemma~4]{Baez2026}}]
\label{prop:walk-model}
Fix $r\in V_G$. A vertex of $G^{\ab}$ may be written $\hml{P}$, where
$P$ is a walk starting at $r$ and ending at $\tau(P)$, with equivalence relation
\begin{equation}\label{eq:equivalence}
  P\equiv Q\quad\Longleftrightarrow\quad
  \tau(P)=\tau(Q)\text{ and }c(P)=c(Q).
\end{equation}
An outgoing arc at $\hml{P}$ is specified by an outgoing arc $a$ at
$\tau(P)$, and joins $\hml{P}$ to $\hml{Pa}$. The projection
$p(\hml{P})=\tau(P)$ is a connected regular cover with
\begin{equation}\label{eq:deck}
  D:=\Deck(p)\cong H_1(G;\Z)\cong\Z^{\beta(G)}.
\end{equation}
A closed walk in $G$ lifts closed if and only if its signed edge-count
vector is zero. The graph $G^{\ab}$ has no loops or parallel edges,
even when $G$ has them.
\end{proposition}

In the maximal abelian cover $p:G^{\ab}\to G$, a nonempty
connected subgraph admits a lifted copy if and only if
it is a tree. Indeed, a cycle has a nonzero signed
edge-count vector and therefore cannot lift to a closed
walk, while trees admit lifted copies by the lifting
criterion above.

The lifting construction for automorphisms applies unchanged to
isomorphisms between two graphs, with the roots chosen compatibly.

\begin{lemma}[{\cite[Lemma~9]{Baez2026}}]\label{lem:lift}
If $f:G\to H$ is an isomorphism, take roots $r$ and $f(r)$. The map
\[
  \hml{P}\longmapsto\hml{f(P)},
\]
with the corresponding map on arcs, is a well-defined isomorphism
$G^{\ab}\to H^{\ab}$.
\end{lemma}

This gives the easy direction of the isomorphism assertions below.

\section{The harmonic realization and its fibres}
\label{sec:coordinates}

Assume for this section that $\beta(G)=b\ge1$, put $X:=G^{\ab}$, fix a root $r\in V_G$ and a lift $x_r\in p^{-1}(r)$ and choose the
fundamental-cycle basis from Lemma~\ref{lem:tree-basis}, denoted by $z_j$ as in the lemma. 
By Proposition~\ref{prop:walk-model}, each vertex of $X$ is represented
by a walk $P$ in $G$ starting at $r$, whose lift from $x_r$
ends at that vertex. Two such walks represent the same vertex
precisely when they have the same endpoint and the same signed
edge-count vector. We may therefore define
\begin{equation}\label{eq:coordinates}
  F_j(\hml{P}):=\langle c(P),z_j\rangle,\qquad
  F:=(F_1,\ldots,F_b):V_X\longrightarrow\Z^b.
\end{equation}
In particular, $F(x_r)=0$.

Note that for a fixed cycle basis, changing the basepoint only adds a constant
vector to $F$, so every difference $F(y)-F(x)$ is independent of this
choice. We fix a basepoint and suppress it from the notation.
Up to an invertible linear change of coordinates, $F$ is the classical
cycle-space realization; see
\cite[Section~4, equations~(28)--(29) and~(35)]{KotaniSunada2001}
and \cite[Section~3, Definition~1]{Baez2026}.\footnote{Let $\pi:\R^{E_G}\to\ker\partial$ be the orthogonal
projection onto the real cycle space. With the same basepoint,
Baez's realization is $i(\hml{P})=\pi c(P)$, and $F=L\circ i$, where
$L(v)=(\langle v,z_1\rangle,\ldots,\langle v,z_b\rangle)$.
The map $L:\ker\partial\to\R^b$ is a linear isomorphism because the
$z_j$ form a basis. For Kotani--Sunada, take the maximal abelian
cover and unit weights. Their circuit form associated with $z_j$
is $\omega_j(a)=\langle\epsilon(a),z_j\rangle$, and evaluating their
path-integral realization on this form gives
$\int_P\omega_j=\langle c(P),z_j\rangle=F_j(\hml{P})$.}

The following properties are standard consequences of the classical
harmonic realization. Kotani and Sunada explicitly define the
circuit forms $\omega_c$ associated with oriented cycles
\cite[p.~12, equation~(35)]{KotaniSunada2001}. They define the realization by
integrating harmonic forms along paths and explicitly verify its
harmonicity \cite[p.~11, equations~(28)--(29) and~(32)]{KotaniSunada2001}.
In our coordinates, the forms take values in $\{-1,0,1\}$, so the
distance bound follows by summing the increments along a shortest
path. We include the elementary verification for completeness.

\begin{lemma}\label{lem:harmonic}
For every pair of vertices $x,y\in V_X$ and every $j$, $|F_j(y)-F_j(x)|\leq d_X(x,y)$.

Each $F_j$ is harmonic on vertices, meaning
\begin{equation}\label{eq:averaging}
  \deg(x)F_j(x)=\sum_{a\in A^+_X(x)}F_j(\tau(a)).
\end{equation}
\end{lemma}

\begin{proof}
Let 
\[
\omega_j:A(G)\longrightarrow\Z,\qquad
\omega_j(a)=\langle\epsilon(a),z_j\rangle.
\]
If $p(\widetilde a)=a$, then
\[
  F_j(\tau(\widetilde a))-F_j(\iota(\widetilde a))=\omega_j(a).
\]
Since $z_j$ is a cycle, $\partial z_j=0$. Write
$z_j=\sum_{e\in E^+}k_e e$, where $E^+$ contains the chosen
orientation of each edge of $G$ and $k_e\in\{-1,0,1\}$ accordingly. Then
\[
  \partial e=\one_{\tau(e)}-\one_{\iota(e)},\qquad
  (\partial z_j)_v=\sum_{\tau(e)=v}k_e-\sum_{\iota(e)=v}k_e.
\]
For an arc $a$ of the edge $e$, the value of $\omega_j(a)$ is
$k_e$ when $a$ agrees with the chosen orientation and $-k_e$
otherwise. Hence
\[
  \sum_{\iota(a)=v}\omega_j(a)
  =\sum_{\iota(e)=v}k_e-\sum_{\tau(e)=v}k_e
  =-(\partial z_j)_v=0.
\]
Using the bijection on outgoing arcs induced by $p$, we obtain
\[
  \sum_{\iota(\widetilde a)=x}
    \bigl(F_j(\tau(\widetilde a))-F_j(x)\bigr)
  =\sum_{\iota(a)=p(x)}\omega_j(a)=0,
\]
which is the harmonic condition. Since
$|\omega_j(a)|=|k_e|\leq1$, summing increments along a shortest path gives the distance bound.
\end{proof}

Let $B$ be the spanning subgraph of $G$ containing exactly its bridges;
isolated vertices are included. Every component of $B$ is a finite
tree. Write $\widetilde B:=p^{-1}(B)$ for the full preimage
of this forest. A \emph{lifted bridge} is a lifted copy of an
edge that is a bridge of $G$; it need not be a bridge of $X$.

The next statement is a covering-space analogue of the
bridge-contraction property of the finite graph Abel--Jacobi
map \cite[Lemma~4.6 and Corollary~4.7]{BakerNorine2007}.
Its bridgeless case is \cite[Theorem~3]{Baez2026}.
The integral-potential argument below is closely related to
Speyer's proof \cite{CycleOrthogonality2016}.

\begin{lemma}\label{lem:fibres}
For $x,y\in V_X$,
\begin{equation}\label{eq:fibres}
  F(x)=F(y)\quad\Longleftrightarrow\quad
  x\text{ and }y\text{ lie in the same component of }\widetilde B.
\end{equation}
Each component of $\widetilde B$ projects isomorphically onto a
component of $B$. In particular, every nonempty fibre of $F$ is a
finite tree.
\end{lemma}

\begin{proof}
($\Leftarrow$) Every cycle vector $z_j$ vanishes on every bridge, since no cycle contains a bridge. Thus all $\omega_j$ are zero on bridges and all $F_j$ are
constant along lifted bridges.

($\Rightarrow$) Write $x=\hml{P}$, $y=\hml{Q}$ and suppose $F(x)=F(y)$. Set
\[
  u:=p(x),\qquad v:=p(y),\qquad c=c(P^{-1}Q)=c(Q)-c(P)\in\Z^{E_G}.
\]
Equality of the coordinates says that $c$ is orthogonal to every real
cycle, since the $z_j$ are also a real basis of $\ker\partial$.
Elementary linear algebra gives
$(\ker\partial)^\perp=\im\partial^{\top}$, so there is a real vertex
potential $h$ such that
\begin{equation}\label{eq:gradient}
  c=\partial^{\top}h,\qquad \partial c=\one_v-\one_u.
\end{equation}
In particular, on an oriented edge $e=s\to t$ we have
$c_e=h(t)-h(s)$, and these differences are integers.

We show that every edge with $c_e\neq0$ is a bridge. Choose a real
number $\lambda$ strictly between its endpoint potentials and unequal
to all vertex potentials, and let
$U=\{w\in V_G:h(w)>\lambda\}$. The coefficient of an edge
crossing this cut is positive when it points into $U$ and negative
when it points out. In the sum of $(\partial c)_w$ over $w\in U$,
an edge inside $U$ contributes zero, while an edge outside $U$
contributes nothing. An edge entering $U$ contributes $c_f$, and
one leaving $U$ contributes $-c_f$. Both contributions are $|c_f|$.
Summing the second identity in
\eqref{eq:gradient} over $U$ therefore gives
\begin{equation}\label{eq:unit-cut}
  0<\sum_{f\text{ crossing }(U,U^c)}|c_f|
     =\sum_{w\in U}(\partial  c)_w=\sum_{w\in U}(\one_v-\one_u)_{w}\leq1.
\end{equation}
All summands are positive integers. Hence the cut has exactly one
edge, namely the chosen edge, which is consequently a bridge.
It follows that $c$ is supported on $B$.

If $u\neq v$, they must lie in the same component of $B$.
Indeed, otherwise summing $\partial c=\one_v-\one_u$ over the
component containing $u$ would give $-1$, whereas the boundary of the restriction of $c$ to that component
has sum zero. Let $R$ be the unique
tree path in that component from $u$ to $v$. The chain $c-c(R)$ is a cycle supported on a forest, so it is
zero. If $u=v$, the same forest argument gives $c=0$, and take $R$
to be empty. Thus $c(Q)=c(PR)$, and \eqref{eq:equivalence} shows
that the lift of $R$ from $x$ ends at $y$. This proves the reverse
implication.

Finally, a connected cover of a tree is a copy of that tree.  Apply this to the
restriction over each component of $B$.
\end{proof}

In the next two lemmas, we analyze the structure of the contracted graph and the faithfulness of automorphisms on graphs without leaves.

\begin{lemma}
\label{lem:contracted-simple}
Contract each component of $\widetilde B$ to a vertex, retaining all
other edges as distinct edge objects. The resulting graph has no
loops or parallel edges.
\end{lemma}

\begin{proof}
A nonbridge edge of a finite connected graph lies on a cycle: join
its endpoints after deleting it, or use the edge itself if it is a
loop. Its signed unit vector therefore has nonzero inner product
with some cycle. By \eqref{eq:coordinates}, a lift of this edge cannot
join two vertices with equal $F$-coordinates. Lemma~\ref{lem:fibres}
then excludes loops in the contracted graph.

Suppose lifted arcs $e:x\to y$ and $f:x'\to y'$ become parallel arcs
between the same ordered pair of contracted vertices. Join $x$ to
$x'$ by a lifted bridge path $R$, and $y'$ to $y$ by another such
path $S$. Then $RfS\bar e$ is a closed walk in $X$. Its projection
has zero signed edge counts by Proposition~\ref{prop:walk-model}.
The nonbridge part of that count is
\[
  \epsilon(p(f))-\epsilon(p(e))=0.
\]
Hence $p(f)=p(e)$ as arcs. Also $x,x'$ are in the same lifted bridge
tree and project to the same vertex. Lemma~\ref{lem:fibres} gives
$x=x'$, and uniqueness of an outgoing lift now gives $e=f$.
\end{proof}

\begin{lemma}
\label{lem:faithful}
Assume $\deg_G(v)\geq2$ for every $v$. If $\ba\in\Aut(X)$ satisfies
$F(\ba x)=F(x)$ for all vertices, then $\ba$ is the identity.
\end{lemma}

\begin{proof}
By Lemma~\ref{lem:fibres}, $\ba$ preserves each component of
$\widetilde B$. An edge of $X$ belongs to $\widetilde B$ exactly
when its endpoints have equal coordinates, by the preceding
proofs. Thus $\ba$ induces an automorphism of the contracted
graph fixing every vertex. Lemma~\ref{lem:contracted-simple}
then fixes every arc of that graph. Consequently every arc
outside $\widetilde B$, and both its endpoints, is fixed.
Since every vertex of $X$ has degree at least two, each
single-vertex component of $\widetilde B$, and each leaf of a
larger component, is incident with an edge outside
$\widetilde B$. These vertices are therefore fixed.

An automorphism of a finite tree fixing every leaf fixes the entire
tree: each vertex lies on a path between leaves, and such a path is
fixed pointwise by the distances from its fixed endpoints. Hence
all vertices of $X$ are fixed. The simplicity of $X$, from
Proposition~\ref{prop:walk-model}, then fixes all arcs as well.
\end{proof}

\section{Recovering the canonical deck group}
\label{sec:deck}

For a connected locally finite graph $X$, define its group of
\emph{bounded-displacement automorphisms}, classically called
\emph{bounded automorphisms} \cite{Trofimov1984,GodsilEtAl1989}, by
\begin{equation}\label{eq:bd}
  \BD(X):=\left\{\ba\in\Aut(X):
       \sup_{x\in V_X}d_X(x,\ba x)<\infty\right\}.
\end{equation}
It is a subgroup: displacement bounds add under composition and
are unchanged by inversion. If $\varphi:X\to X'$ is an isomorphism,
preservation of graph distance gives
\begin{equation}\label{eq:bd-intrinsic}
  \varphi\BD(X)\varphi^{-1}=\BD(X').
\end{equation}
Thus this group is intrinsic to the abstract graph, independently
of a choice of covering projection.

We use the ordinary maximum principle for harmonic functions
\cite[Section 2, p.~5]{KotaniSunada2001} in the following elementary form.

\begin{lemma}\label{lem:maximum}
Let $u$ be a real-valued function on the vertices of a connected graph $X$
of positive finite degrees. Suppose its image is finite and
\[
  \deg(x)u(x)=\sum_{a\in A_X^+(x)}u(\tau(a))
\]
at every vertex. Then $u$ is constant.
\end{lemma}

\begin{proof}
The finite image has a maximum $M$, attained at some vertex. All
neighbouring values are at most $M$, so their average can equal $M$
only if they all equal $M$. Connectivity propagates this conclusion
to every vertex.
\end{proof}

For finite $c$-connected base graphs of cycle rank at
least two, the following identification is established
in Eon's proof of Theorem~6.1
\cite[pp.~62--63]{Eon2007}.
We give a self-contained proof under the weaker
hypothesis of minimum degree at least two, using
integral harmonic coordinates and a comparison of
cycle ranks.

\begin{theorem}
\label{thm:deck}
Let $G$ be finite and connected with minimum degree at least two and
$b=\beta(G)\geq2$. For its maximal abelian cover $p:X\to G$,
\begin{equation}\label{eq:deck-recovery}
  \BD(X)=\Deck(p).
\end{equation}
\end{theorem}

\begin{proof}
Put $\Gamma:=\BD(X)$ and $D:=\Deck(p)\cong\Z^b$.
Use the Kotani--Sunada realization \cite[Section 4]{KotaniSunada2001}
in the integral coordinates $F$
of \eqref{eq:coordinates}.
The relevant properties were proved in Section~\ref{sec:coordinates}:
the coordinates have bounded edge increments and are harmonic,
and coordinate-preserving automorphisms are trivial when there are
no leaves. We now use these facts to recover $D$ from $X$.

\medskip\noindent\textit{Step 1: bounded displacement gives a constant
coordinate translation.}
For $\ba\in \Gamma$, choose an integer $R$ such that $d_X(x,\ba x)\leq R$ for
every $x$. For each $j$ set $u_j(x)=F_j(\ba x)-F_j(x)$. An automorphism
preserves degrees and permutes outgoing arcs, so both $F_j\circ \gamma$
and $F_j$ satisfy the averaging equation in
Lemma~\ref{lem:harmonic}. Hence $u_j$ does too. The same lemma gives
\[
  u_j(x)\in\{-R,-R+1,\ldots,R\}.
\]
Lemma~\ref{lem:maximum} makes $u_j$ constant. Therefore there is a
vector $t_\ba\in\Z^b$ such that
\begin{equation}\label{eq:translation}
  F(\ba x)=F(x)+t_\ba\qquad(x\in V_X).
\end{equation}
Composition gives $t_{\ba \ba'}=t_\ba+t_{\ba'}$. By
Lemma~\ref{lem:faithful}, $t_\ba=0$ implies $\ba=1$. Thus
\begin{equation}\label{eq:translation-injection}
\begin{split}
    t:\Gamma\longrightarrow\Z^b\\
    \ba\longmapsto t_\ba
\end{split}
\end{equation}
is an injective homomorphism. In particular, $\Gamma$ is free abelian. If $\ba$ fixes a vertex, \eqref{eq:translation} gives
$t_\ba=0$, so $\ba=1$. If $\ba$ inverts an edge, then $\ba^2$ fixes its
endpoints, so $\ba^2=1$. Since $\Gamma$ is torsion-free, $\ba=1$, a contradiction.
Hence $\Gamma$ acts freely on vertices and without edge inversions.

\medskip\noindent\textit{Step 2: the canonical deck group is contained
in $\Gamma$ and has finite index.}
Choose a lift $x_v$ above each of the finitely many base vertices.
For $d,d'\in D$, commutativity and preservation of distance give
\[
  d_X(d'x_v,dd'x_v)=d_X(x_v,dx_v).
\]
The $D$-orbits of the $x_v$ exhaust $V_X$. Hence the displacement
of $d$ is at most $\max_v d_X(x_v,dx_v)$, and $D\leq \Gamma$.

Since $\Gamma$ is abelian, $D$ is normal in $\Gamma$. Fix a vertex $x\in X$.
The map
\[
  \Gamma/D\longrightarrow V_G,\qquad \ba D\longmapsto p(\ba x)
\]
is well defined and injective. If $\ba D=\ba 'D$, then $\ba=d\ba'$ for some
$d\in D$, and $p(\ba x)=p(d\ba'x)=p(\ba'x)$. Conversely, if
$p(\ba x)=p(\ba'x)$, regularity gives $d\in D$ with $\ba x=d\ba'x$. Then $(d\ba')^{-1}\ba$ fixes $x$, so it is the
identity by Step~1. Hence $\ba D=\ba'D$. It follows that
\begin{equation}\label{eq:index}
  k:=[\Gamma:D]\leq |V_G|<\infty.
\end{equation}

\medskip\noindent\textit{Step 3: form the finite covering quotient.}
Let $Y:=X/\Gamma$. Since $\Gamma$ acts freely on vertices and without
edge inversions, the projection is a regular covering with deck group
$\Gamma$; see \cite[Theorem~2.2.2]{GrossTucker1987}.
As $D\trianglelefteq \Gamma$, this projection factors through $G=X/D$:
\[
  X\longrightarrow G\longrightarrow Y.
\]

The induced map $G\to Y$ is a covering, since the local arc maps from $X$ to
both quotients are bijections. Its vertex fibres have size $k$,
being the $D$-orbits within a free $\Gamma$-orbit. The graph $Y$ is finite
and connected. Counting vertices and outgoing arcs in the
$k$-sheeted cover gives:
\[
  |V_G|=k|V_Y|,\qquad |E_G|=k|E_Y|,
\]
also when loops or parallel edges occur. Consequently
\begin{equation}\label{eq:rank-formula}
  b-1=k\bigl(\beta(Y)-1\bigr).
\end{equation}

\medskip\noindent\textit{Step 4: the abelian deck group forces enough
cycles in $Y$.}
Fix $x\in X$ above $y\in Y$. A closed walk at $y$ has a lift from
$x$ ending at $\ba x$ for a unique $\ba\in \Gamma$. This defines a surjective
homomorphism $\pi_1(Y,y)\to \Gamma$. It is well defined because backtracks
lift to backtracks; concatenation gives multiplication because the
second lift is obtained by translating a lift from $x$. It is onto
because $X$ is connected, so a path from $x$ to any $\ba x$ exists.
As $\Gamma$ is abelian, Lemma~\ref{lem:tree-basis} gives a surjection
\begin{equation}\label{eq:homology-surjection}
  H_1(Y;\Z)\cong\Z^{\beta(Y)}\longtwoheadrightarrow{} \Gamma.
\end{equation}

The inclusion $D\hookrightarrow \Gamma$ and the map $t$ give
\[
  D\cong\Z^b{\lhook\joinrel\longrightarrow} \Gamma\overset{t}{\lhook\joinrel\longrightarrow} \Z^b,
  \qquad b=\rank_\Z D\leq\rank_\Z \Gamma\leq b.
\]
The surjection \eqref{eq:homology-surjection} implies
\begin{equation}\label{eq:rank-lower-bound}
  \begin{split}
  \rank_\Z \Z^{\beta(Y)} &\geq \rank_\Z \Gamma \\    
  \beta(Y)&\geq b
  \end{split}
\end{equation}
Combining \eqref{eq:rank-formula} and
\eqref{eq:rank-lower-bound},
\[
  b-1=k\bigl(\beta(Y)-1\bigr)\geq k(b-1).
\]
Since $b-1>0$ and $k$ is a positive integer, $k=1$. Hence $\Gamma=D$.
\end{proof}

We observe that the differences $F_j\circ \ba-F_j$ are bounded; the coordinates
$F_j$ themselves are generally unbounded. A bounded real-valued
function need not attain its supremum on an infinite graph. Here the
differences are integer-valued, so boundedness makes their image
finite and the elementary maximum principle applies.

The following descent statement is an immediate consequence
of Theorem~\ref{thm:deck}.

\begin{corollary}\label{cor:core-rigidity}
Let $G,H$ be finite connected graphs of minimum degree at least two
and cycle rank at least two. Every isomorphism
$\varphi:G^{\ab}\to H^{\ab}$ descends to a unique isomorphism
$f:G\to H$ satisfying
\[
  p_H\circ\varphi=f\circ p_G.
\]
\[
\begin{tikzcd}[row sep=large, column sep=large]
G^{\mathrm{ab}} \arrow[r, "\varphi"] \arrow[d, "p_G"]
  & H^{\mathrm{ab}} \arrow[d, "p_H"] \\
G \arrow[r, "f"] & H
\end{tikzcd}
\]
\end{corollary}

\begin{proof}
Equation~\eqref{eq:bd-intrinsic} and Theorem~\ref{thm:deck} give
\[
  \varphi\Deck(p_G)\varphi^{-1}=\Deck(p_H).
\]
For $v\in V_G$, choose $x\in p_G^{-1}(v)$ and define
$f(v)=p_H(\varphi(x))$. If $x'$ is another lift of $v$, regularity
gives $x'=gx$ for some $g\in\Deck(p_G)$. Then
$h:=\varphi g\varphi^{-1}\in\Deck(p_H)$, so
\[
  p_H(\varphi(x'))=p_H(h\varphi(x))=p_H(\varphi(x)).
\]
Thus $f$ is well defined. The same construction on arc orbits
preserves incidence and reversal. Applying it to $\varphi^{-1}$
gives the inverse graph map, so $f$ is an isomorphism.
Finally, $p_H\circ\varphi=f\circ p_G$, and the surjectivity of
$p_G$ makes $f$ unique.
\end{proof}

An automorphism of a cover \emph{descends} when it sends each
covering fibre onto a covering fibre and hence induces an automorphism
of the base. Following \cite[Section~5]{Baez2026}, let $\Cov(G)$ be
the subgroup of $\Aut(G^{\ab})$ formed by these automorphisms.
The next corollary says that, for a base without leaves, every
symmetry of the cover has this property.
\begin{corollary}\label{cor:all-descend}
Let $G$ be finite and connected with minimum degree at least two.
For $p:G^{\ab}\to G$,
\[
  \Cov(G)=\Aut(G^{\ab}).
\]
There is a short exact sequence
\begin{equation}\label{eq:full-aut-sequence}
  1\longrightarrow H_1(G;\Z)\overset{\iota}{\longrightarrow}\Aut(G^{\ab})
  \overset{\rho}{\longrightarrow}\Aut(G)\longrightarrow1,
\end{equation}
The injection $\iota$ is the one used by Baez
\cite[Section~5]{Baez2026} to identify $H_1(G;\Z)$ with the deck group.
In the walk model rooted at $r$, it is given by
\[
  \iota(z)(\hml{P})=\hml{C_zP},
  \qquad z\in H_1(G;\Z),
\]
where $C_z$ is any closed walk based at $r$ with $c(C_z)=z$.
On arcs, it sends the arc specified by $(\hml{P},a)$ to the one
specified by $(\hml{C_zP},a)$. Thus $\iota(z)$ is the deck
transformation that adds $z$ to the signed edge-count vector.
The map $\rho$ sends each automorphism to the one it induces on $G$.
Explicitly, every base automorphism lifts, and an automorphism
upstairs induces the identity downstairs exactly when it is a deck
transformation. Thus any two lifts of the same base automorphism
differ by a deck transformation.
\end{corollary}

\begin{proof}
If $\beta(G)\geq2$, apply Corollary~\ref{cor:core-rigidity} with $H=G$.
Otherwise, handshake lemma gives $\beta(G)\geq1$, and equality
forces every degree to be two. Thus $G$ is a cycle of some length $m$.
Its maximal abelian cover is the two-way infinite path, with vertices
labelled by $\Z$ and projection modulo $m$. Every automorphism has the
form $n\mapsto\varepsilon n+c$, where $\varepsilon\in\{1,-1\}$ and
$c\in\Z$, and therefore descends, including on arcs when $m=1$ or $2$.

The identity $c(C_zP)=z+c(P)$ shows that $\iota$ is well defined
and that $\iota(z+w)=\iota(z)\circ\iota(w)$. Its action on the
class of the empty walk at $r$ shows that $\iota(z)$ can be the
identity only when $z=0$, so $\iota$ is injective. Every deck
transformation is determined by the image of this vertex, and every
vertex above $r$ is represented by a closed walk. Hence
$\iota(H_1(G;\Z))=\Deck(p)$.
The descended map $\rho$ is unique and respects composition, with
\[
  \ker\rho=\Deck(p)=\iota(H_1(G;\Z)).
\]
To prove surjectivity, let $f\in\Aut(G)$ and choose a walk $Q$ from
$r$ to $f(r)$. The map $\hml{P}\mapsto\hml{Qf(P)}$ is a well-defined
automorphism of the walk model covering $f$; reversing $Q$ and applying
$f^{-1}$ constructs its inverse. Hence $\rho$ is onto.
\end{proof}

For $\Cov(G)$ the exact sequence is already
\cite[Lemma 9]{Baez2026}. The conclusion here identifies that group
with the full abstract automorphism group. Splitting of the sequence
is a separate question. For the classical lifting problem along a
specified regular cover, and the structure of lifted groups, see
\cite[Sections 7 and 9]{MalnicEtAl2000}.

\section{Attached trees and the small cycle ranks}
\label{sec:classification-proof}

In this section, we analyze the decomposition of a finite
connected graph into its core and attached trees, and treat
the cases of small cycle rank, completing the proof of
Theorem~\ref{thm:classification}.

Recall that the \emph{core}, or \emph{$2$-core}, of a finite
graph $G$, denoted by $\Core(G)$, is obtained by repeatedly
deleting all vertices of degree less than two, together with
their incident edges, until no such vertices remain. Degrees
are recomputed in the remaining graph at each step. The core is empty if and only if the graph is a forest. Otherwise it has minimum degree at least two and is connected when the graph is.

\begin{lemma}\label{lem:core}
Let $G$ be finite and connected with $\beta(G)\geq1$ and
let $C:=\Core(G)$. Then:
\begin{enumerate}[(i)]
\item $C$ is finite, connected, has minimum degree at least two, and
      $\beta(C)=\beta(G)$.
\item $G$ is obtained from $C$ by attaching a finite rooted tree
      $T_v$ at each $v\in V_C$, identifying only its root with $v$.
      A tree consisting of just the root is allowed.
\item In $X=G^{\ab}$, the inverse image $p^{-1}(C)$ is connected and
      canonically a copy of $C^{\ab}$ as a cover of $C$. At each
      vertex above $v$, a copy of $T_v$ is attached.
\item Applying the same iterative deletion to $X$ leaves precisely
      $p^{-1}(C)$. Thus every graph isomorphism between two such
      covers preserves the preimages of their cores and their
      rooted decorations.
\end{enumerate}
\end{lemma}

\begin{proof}
Since $G$ has a cycle, the deletion process never removes all its
cycles. Each individual leaf deletion removes one vertex and one
edge, preserves connectedness, and leaves $\beta$ unchanged. There
is no isolated vertex in the remaining connected graph containing
a cycle. The finite process therefore ends at a nonempty connected
graph of minimum degree at least two. Reversing the deletions
attaches one new leaf at a time; this proves (i) and (ii).

Take the root of the walk model in $C$. Every walk from $C$
into an attached tree returns to the same core vertex and reduces
by backtracks, so it has zero chain. A walk ending in $C$ therefore
represents the same covering vertex as a walk entirely in $C$.
Moreover, two walks in $C$ have equal chains in $G$ exactly when
they have equal chains in $C$. The walk model hence identifies
$p^{-1}(C)$ with $C^{\ab}$ and proves its connectedness. The
restriction of a cover to a tree consists of copies of that tree,
as proved in Lemma~\ref{lem:fibres}; this gives the decorations in
(iii).

For (iv), let $G_n$ and $X_n$ denote the graphs after $n$ simultaneous
rounds of deletion, starting with $G$ and $X$. Inductively,
$X_n=p^{-1}(G_n)$, and the restricted map remains a covering on the
remaining outgoing arcs. Thus a vertex of $X_n$ has degree less
than two precisely when its image in $G_n$ does. Once $G_n=C$,
the process is stable upstairs as well. The procedure is invariant
under graph isomorphisms. The attached rooted trees are then
recovered as the root together with the components outside the
core incident with it.
\end{proof}

\begin{proof}[Proof of Theorem~\ref{thm:classification}, case (c)]
Assume $\beta(G),\beta(H)\geq2$ and let $\varphi$ be an isomorphism
of their maximal abelian covers $X=G^{\ab}$ and $Y=H^{\ab}$.
Put $C_G=\Core(G)$ and
$C_H=\Core(H)$.
Lemma~\ref{lem:core} restricts $\varphi$ to
$C_G^{\ab}\cong C_H^{\ab}$. By
Corollary~\ref{cor:core-rigidity}, this restriction descends to an
isomorphism $f:C_G\to C_H$.

Write $\widetilde C_G=p_G^{-1}(C_G)$ and
$\widetilde C_H=p_H^{-1}(C_H)$ for the preimages of the cores.
For $v\in V_{C_G}$, choose $x\in p_G^{-1}(v)$.
Since the restriction of $\varphi$ to these preimages descends
to $f$, we have $p_H(\varphi(x))=f(v)$.

The rooted tree attached at $x$ consists of $x$, together with
all components of $X-V_{\widetilde C_G}$ adjacent to $x$ and
their attaching edges. Each such component is attached to the
preimage of the core at exactly one vertex.
Because $\varphi(\widetilde C_G)=\widetilde C_H$, the isomorphism
$\varphi$ maps components of $X-V_{\widetilde C_G}$ bijectively
onto components of $Y-V_{\widetilde C_H}$, preserving their
attachment vertices. Consequently, it maps the components
attached at $x$ precisely onto those attached at $\varphi(x)$.
Thus $\varphi$ restricts to an isomorphism of the corresponding
rooted trees.

Thus $T_v\cong T_{f(v)}$ as rooted trees. Choose such an isomorphism for
each core vertex and combine them with $f$. Since the attached
trees meet the core only in their roots, these maps give
$G\cong H$. The converse is Lemma~\ref{lem:lift}.
\end{proof}

The preceding argument reconstructs the base and its attached trees,
but an individual isomorphism of the covers need not descend.
For example, if two leaves are attached at a core vertex, swapping
them at just one vertex upstairs gives a bounded-displacement
automorphism that neither descends nor belongs to the deck group.
This motivates the first problem in Section~\ref{sec:baez}:
Corollary~\ref{cor:abelian-cover-descent} gives a complete
characterization of when every automorphism descends. Apart from
the tree case, this happens exactly when each attached tree has
trivial root-preserving automorphism group. Section~\ref{sec:baez-splitting}
then addresses the second problem, giving in
Theorem~\ref{thm:baez-orbit-splitting} a criterion for choosing the
lifts compatibly with composition.

\begin{proof}[Completion of the proof of Theorem~\ref{thm:classification}]
If $\beta(G)=0$, the connected graph $G$ is a finite tree. In the
walk model every walk is equivalent to the unique simple path to
its endpoint, so $G^{\ab}=G$. This proves (a).

If $\beta(G)=1$, Lemma~\ref{lem:core} gives a finite connected core
$C$ with minimum degree at least two and $|E_C|=|V_C|$. The sum
of the degrees is $2|V_C|$, so every degree is exactly two. Thus
$C$ is a cycle, including the loop and double-edge possibilities.
Its maximal abelian cover is the two-way infinite path: the
fundamental group is $\Z$, its commutator subgroup is trivial,
and the cycle is unwrapped by the covering map. Consequently
$G^{\ab}$ is this path, indexed by $\Z$, with the tree
$T_{i\bmod m}$ attached at position $i$.

The infinite path is the preimage of the core and is intrinsically
determined by Lemma~\ref{lem:core}. Every isomorphism between two such covers
therefore induces an isomorphism of their paths. An isomorphism
of the two-way infinite path has the form
\[
  i\longmapsto \varepsilon i+k,\qquad
  \varepsilon\in\{1,-1\},\quad k\in\Z:
\]
the images of positions $0$ and $1$ determine the rest successively.
It extends to the decorated paths exactly when all corresponding
rooted trees are isomorphic.

To translate this into finite words, the periods of a bi-infinite
decoration sequence form a subgroup of $\Z$ containing $m\Z$.
Thus they equal $\ell\Z$ for a least positive period $\ell$ dividing
$m$. The block of length $\ell$ is the shortest repeating block.
Two sequences agree after a shift or reflection exactly when these
blocks agree up to cyclic rotation or reversal. Conversely, such
agreement gives a path isomorphism, and rooted-tree isomorphisms
at each position extend it to the whole cover. This proves (b).

Finally, the three cases cannot mix. A rank-zero base has a finite
cover equal to itself. A rank-one base has an infinite tree as its
maximal abelian cover. If $\beta(G)\geq2$, take two distinct free
generators $s,t$ of $\pi_1(G,r)$. Their commutator
$sts^{-1}t^{-1}$ is nontrivial and has zero edge-count vector.
Choose its nonempty reduced closed-walk representative. It lifts
to a closed walk in $G^{\ab}$, and that lift is still reduced,
since an immediate backtrack upstairs would project to one
downstairs. A tree has no nonempty reduced closed walk, so this
cover contains a cycle. It is infinite by its deck group
$\Z^{\beta(G)}$. These intrinsic distinctions separate the cases.
Within the last case, the already proved implication $G\cong H$
also forces equality of cycle ranks.
\end{proof}

The final inequality in Theorem~\ref{thm:deck} uses $b-1>0$.
For a cycle of length $m$, the cover is the two-way infinite path;
its bounded-displacement automorphisms are all translations,
whereas the canonical deck group consists of translations by
multiples of $m$. The index is $m$, and the rank formula reads
$0=m\cdot0$. This explains the exception within the proof itself.

\section{Baez's problems}\label{sec:baez}
The tools developed in this article also allow us to answer two problems Baez proposed \cite[Section~8]{Baez2026} about the automorphism group of maximal abelian covers. We present their formulation and solution here.
\subsection{When does every automorphism descend?}
Baez asks when every automorphism of a maximal abelian cover
descends to an automorphism of its base.
The descent theorem for graphs of minimum degree at least two, together with the core decomposition, gives the following criterion for finite connected graphs.

\begin{corollary}
\label{cor:abelian-cover-descent}
Let \(G\) be a finite connected graph, and let
\[
    p \colon X=G^{\mathrm{ab}} \longrightarrow G
\]
be its maximal abelian cover. Define
\[
    \operatorname{Cov}(G)
    :=
    \left\{
        \varphi\in\operatorname{Aut}(X)
        \;\middle|\;
        \begin{array}{l}
        \text{there exists }f\in\operatorname{Aut}(G)\\
        \text{such that }p\circ\varphi=f\circ p
        \end{array}
    \right\}.
\]

Equivalently, the following diagram commutes for some
$f\in\operatorname{Aut}(G)$:
\[
\begin{tikzcd}[row sep=large, column sep=large]
X \arrow[r, "\varphi"] \arrow[d, "p"']
  & X \arrow[d, "p"] \\
G \arrow[r, "f"']
  & G
\end{tikzcd}
\]

Let \(C:=\Core(G)\). When \(C\neq\varnothing\), let \(T_v\) be the component
of \(G-E(C)\) containing \(v\), rooted at \(v\), for each
\(v\in V(C)\). Thus \(G\) is obtained from \(C\) by attaching
the finite rooted trees \(T_v\).

Then
\[
    \operatorname{Cov}(G)=\operatorname{Aut}(X)
\]
if and only if either
\begin{enumerate}
    \item \(G\) is a tree; or
    \item \(C\neq\varnothing\) and
    \[
        \operatorname{Aut}(T_v;v)=\{1\}
        \qquad\text{for every }v\in V(C),
    \]
    where \(\operatorname{Aut}(T_v;v)\) denotes the group
    of automorphisms of \(T_v\) fixing its root.
\end{enumerate}
Whenever an automorphism of \(X\) descends, its descended
automorphism is unique.
\end{corollary}

\begin{proof}
If $G$ is a tree, then $p:X\to G$ is an isomorphism,
so every automorphism of $X$ descends. Henceforth assume
that $\beta(G)\geq1$, and put $C=\Core(G)$.

By Lemma~\ref{lem:core}, the preimage of the core
\[
  Y:=p^{-1}(C)
\]
is connected, and
\[
  q:=p|_Y:Y\longrightarrow C
\]
is the maximal abelian cover of $C$. At each vertex
$x\in V(Y)$, a finite rooted tree $\widetilde T_x$ is
attached, and the restriction
\[
  p_x:=p|_{\widetilde T_x}:
  \widetilde T_x\longrightarrow T_{p(x)}
\]
is a rooted graph isomorphism. Indeed, every closed walk in
$T_{p(x)}$ reduces to the empty walk, so its inclusion into $G$
has a unique lift sending its root to $x$
(Section~\ref{subsec:background}). The image of this lift is
the component of $p^{-1}(T_{p(x)})$ containing $x$, namely
$\widetilde T_x$, and the lift is inverse to $p_x$.

By Lemma~\ref{lem:core}(iv), repeatedly deleting vertices of
degree less than two from $X$ stabilizes after finitely many
rounds at $Y$. Thus $Y$ is determined by the abstract graph $X$,
independently of $p$. Every $\varphi\in\operatorname{Aut}(X)$
therefore preserves $Y$ and maps $\widetilde T_x$ onto
$\widetilde T_{\varphi(x)}$.

Suppose first that every attached tree is root-rigid:
\[
  \operatorname{Aut}(T_v;v)=\{1\}
  \qquad(v\in V(C)).
\]
Let $\varphi\in\operatorname{Aut}(X)$. Since $C$ has minimum
degree at least two, Corollary~\ref{cor:all-descend} gives
an automorphism $f_C\in\operatorname{Aut}(C)$ for which
the following diagram commutes:
\begin{equation}\label{eq:baez-core-descent}
\begin{tikzcd}[row sep=large, column sep=large]
Y \arrow[r, "\varphi|_Y"] \arrow[d, "q"']
  & Y \arrow[d, "q"] \\
C \arrow[r, "f_C"']
  & C .
\end{tikzcd}
\end{equation}
We will extend $f_C$ over the attached trees.

For $v\in V(C)$ and any lift $x\in p^{-1}(v)$, define
the rooted isomorphism
\[
  \psi_x
  :=p_{\varphi(x)}
     \circ\bigl(\varphi|_{\widetilde T_x}\bigr)
     \circ p_x^{-1}
  :T_v\longrightarrow T_{f_C(v)}.
\]
The target is as stated because commutativity of
\eqref{eq:baez-core-descent} gives
$p(\varphi(x))=f_C(p(x))=f_C(v)$. By definition of $\psi_x$,
the following diagram also commutes:
\begin{equation}\label{eq:baez-tree-descent}
\begin{tikzcd}[row sep=large, column sep=huge]
\widetilde T_x
  \arrow[r, "\varphi|_{\widetilde T_x}"] \arrow[d, "p_x"']
  & \widetilde T_{\varphi(x)} \arrow[d, "p_{\varphi(x)}"] \\
T_v \arrow[r, "\psi_x"']
  & T_{f_C(v)} .
\end{tikzcd}
\end{equation}

This induced map does not depend on the choice of $x$.
Indeed, if $x'$ is another lift of $v$, then $\psi_x$ and
$\psi_{x'}$ are rooted isomorphisms between the same two
trees. Therefore
\[
  \psi_{x'}^{-1}\circ\psi_x
  \in\operatorname{Aut}(T_v;v)=\{1\},
\]
so $\psi_x=\psi_{x'}$. Denote their common value by $\psi_v$.

Define $f:G\to G$ by
\[
  f|_C=f_C,
  \qquad
  f|_{T_v}=\psi_v
  \quad(v\in V(C)).
\]
These definitions agree at the roots, since
$\psi_v(v)=f_C(v)$. Moreover, $f_C$ permutes the core
vertices, and each $\psi_v$ maps $T_v$ isomorphically onto
$T_{f_C(v)}$. Thus $f$ is a graph automorphism.

Commutativity of \eqref{eq:baez-core-descent} gives
$p\circ\varphi=f\circ p$ on $Y$. On each attached tree
$\widetilde T_x$, commutativity of \eqref{eq:baez-tree-descent}
and $\psi_x=\psi_{p(x)}$ give the same identity.
Hence $p\circ\varphi=f\circ p$ on all of $X$, proving
that every automorphism descends.

Conversely, suppose that some $T_v$ admits a nonidentity
root-preserving automorphism $\alpha$. Choose
$x\in p^{-1}(v)$ and define an automorphism $k$ of $X$
by
\[
  k|_{\widetilde T_x}=p_x^{-1}\circ\alpha\circ p_x,
\]
taking $k$ to be the identity elsewhere. This is
well defined because $\alpha$ fixes the root and
$\widetilde T_x$ meets the rest of $X$ only at $x$.

Since $\beta(G)\geq1$, the deck group
$H_1(G;\mathbb Z)$ is infinite, so there is another
lift $x'\in p^{-1}(v)$ with $x'\neq x$.
Choose a vertex $w\in V(T_v)$ moved by $\alpha$, and put
\[
  z=p_x^{-1}(w),
  \qquad
  z'=p_{x'}^{-1}(w).
\]
These vertices lie in the same covering fibre, but
\[
  p(k(z))=\alpha(w)\neq w=p(k(z')).
\]
Thus $p\circ k$ is not constant on that fibre, so $k$
cannot descend. This proves necessity.

Finally, whenever $p\circ\varphi=f\circ p$, the
surjectivity of $p$ determines $f$ uniquely.
\end{proof}

The criterion allows cut vertices, bridges, and leaves.
In particular, every graph of minimum degree at least two
satisfies it, since each attached rooted tree consists
only of its root.

For graphs with nonempty core, attached trees obstruct descent
precisely when one has a nontrivial root-preserving automorphism.
A path rooted at an endpoint causes no obstruction, whereas
two interchangeable leaves at the same root do.

In the unicyclic case, no assumption that the decoration
word is primitive is needed: every shift or reflection
of the preimage of the core normalizes the subgroup of
translations by the original cycle length.
\subsection{When does the lift extension split?}
\label{sec:baez-splitting}

Every automorphism of a finite connected graph lifts to its maximal
abelian cover. Choosing these lifts independently, however, need not
give a group action: the chosen lift of a product may differ from the
product of the chosen lifts by a deck transformation. Baez's second
question concerns precisely this compatibility problem
\cite[Section~8]{Baez2026}. We give a criterion expressed entirely in
terms of the vertex and edge orbits of the base graph.

Recall that the group of covering automorphisms is
\[
 \operatorname{Cov}(G):
 =\{f\in\operatorname{Aut}(X):p\circ f=h\circ p
       \text{ for some }h\in \operatorname{Aut}(G)\}.
\]
The automorphism $h$ is unique. Sending $f$ to $h$ defines the canonical
exact sequence \cite[Lemma~9]{Baez2026}
\begin{equation}
 1\longrightarrow H_1(G;\mathbb Z)\xlongrightarrow{\iota}\operatorname{Cov}(G)
   \xlongrightarrow{\rho}\operatorname{Aut}(G)\longrightarrow 1.
 \label{eq:baez-lift-extension}
\end{equation}
Here $\iota$ is the inclusion of $H_1(G;\mathbb Z)$ as the deck group
used by Baez \cite[Section~5]{Baez2026}, explicitly described in
Corollary~\ref{cor:all-descend}. The extension
\emph{splits} if there is a group homomorphism
$s\colon\operatorname{Aut}(G)\to\operatorname{Cov}(G)$ such that
$\rho\circ s=\operatorname{id}_{\operatorname{Aut}(G)}$.
Such a homomorphism $s$ is called a \emph{section} of $\rho$.
\[
\begin{tikzcd}[column sep=1.8em]
1 \arrow[r] &
H_1(G;\mathbb Z) \arrow[r, hook, "\iota"] &
\operatorname{Cov}(G) \arrow[r, two heads, "\rho"] &
\operatorname{Aut}(G) \arrow[r] \arrow[l, bend left=45, "s"'] &
1
\end{tikzcd}
\]

Thus a splitting is a compatible choice of lifts of all
elements of $\operatorname{Aut}(G)$.
As usual, the extension splits if and only if it is isomorphic,
as an extension, to the semidirect product
$H_1(G;\mathbb Z)\rtimes\operatorname{Aut}(G)$ with the natural
homology action \cite[Section~4, p.~122]{MalnicPozar2016}.
Here even an abstract group isomorphism suffices; see
Corollary~\ref{cor:baez-semidirect}.

Malni\v{c} and Po\v{z}ar give general splitting criteria,
including an integer-linear-system criterion for finitely
presented abelian deck groups
\cite[Theorems~4.1 and~5.2]{MalnicPozar2016}.
These results include maximal abelian covers; here we obtain
an explicit criterion using only vertex and unoriented-edge
orbit sizes. Bridson and Vogtmann \cite[Section~4.1]{BridsonVogtmann2012}
give examples of finite cyclic groups of graph automorphisms
whose actions do not lift to their maximal abelian covers.

Let $V_1,\ldots,V_q$ be the $\operatorname{Aut}(G)$-orbits on vertices, and let
$E_1,\ldots,E_t$ be the $\operatorname{Aut}(G)$-orbits on \emph{unoriented edges}.
Parallel edges are counted separately, and a loop counts as one
unoriented edge. We write $d_{\mathrm{orb}}(G)$ for the greatest
common divisor of all these orbit sizes:
\begin{equation}
 d_{\mathrm{orb}}(G):=\gcd\bigl(|V_1|,\ldots,|V_q|,
                         |E_1|,\ldots,|E_t|\bigr).
 \label{eq:baez-orbit-gcd}
\end{equation}
If $G$ has no edges, the edge-orbit list is simply omitted.

\begin{theorem}
\label{thm:baez-orbit-splitting}
The canonical extension \eqref{eq:baez-lift-extension} splits if and
only if $d_{\mathrm{orb}}(G)=1$. When it splits, a choice of section gives an
isomorphism
\[
\begin{aligned}
 H_1(G;\mathbb Z)\rtimes\operatorname{Aut}(G)
   &\longrightarrow\operatorname{Cov}(G),\\
 (a,h)&\longmapsto\iota(a)s(h),
\end{aligned}
\]
where $\operatorname{Aut}(G)$ acts on $H_1(G;\mathbb Z)$ by its natural action on integral homology.
This isomorphism identifies the first factor with the canonical deck
group and intertwines projection onto $\operatorname{Aut}(G)$ with $\rho$.
\end{theorem}

We prepare the proof of Theorem~\ref{thm:baez-orbit-splitting}
in three steps. First, we describe the group of lifts using
integral chains. Second, we prove an elementary fact about
integral vectors under coordinate permutations. These two
ingredients show that the extension splits precisely when the
barycentric subdivision admits an invariant integral $0$-chain
of degree one. We then complete the proof by determining the
possible degrees of such invariant chains. The final part of
the proof of Theorem~\ref{thm:baez-orbit-splitting} appears on
page~\pageref{proof:baez-orbit-splitting}.

An automorphism of $G$ may reverse an edge, which introduces signs
in its action on integral $1$-chains. To remove this difficulty, let
$S:=S(G)$ be the barycentric subdivision: insert one new vertex $m_e$
in each unoriented edge $e$. Orient each resulting half-edge from its
original endpoint toward $m_e$. Every element of $\operatorname{Aut}(G)$ preserves this
orientation. A loop produces two distinct parallel half-edges with the
same orientation; an automorphism reversing the loop interchanges
these half-edges.

Each automorphism $h\in\operatorname{Aut}(G)$ induces an
automorphism of $S$, also denoted by $h$: it sends an original
vertex $u$ to $h(u)$, an edge midpoint $m_e$ to $m_{h(e)}$,
and each oriented half-edge to the corresponding oriented
half-edge. We extend these maps $\mathbb Z$-linearly to chains.
Explicitly, for an integral $1$-chain
$c=\sum_a n_a a\in C_1(S;\mathbb Z)$ and an integral $0$-chain
$z=\sum_v b_vv\in C_0(S;\mathbb Z)$, define
\[
  hc:=\sum_a n_a h(a),
  \qquad
  hz:=\sum_v b_v h(v).
\]
Here $a$ ranges over the oriented half-edges chosen above,
$v$ ranges over $V(S)$, and all coefficients are integers.
Thus $hc$ denotes the action of the automorphism $h$ on the
chain $c$. These are left actions: $(hk)c=h(kc)$ and $1c=c$,
where $hk:=h\circ k$.

Consequently, both chain groups have bases permuted by
$\operatorname{Aut}(G)$. The boundary map is
\[
 \partial\colon C_1(S;\mathbb Z)\longrightarrow C_0(S;\mathbb Z),
 \qquad \partial(u\to m_e)=m_e-u.
\]
Because $h$ preserves the initial and terminal vertices of every
oriented half-edge, linearity gives
\[
  \partial(hc)=h(\partial c).
\]

Subdivision preserves the fundamental group and its abelianization.
Moreover, subdividing each edge of $X$ gives the maximal abelian
cover of $S$. A lift of $h$ subdivides uniquely to a lift of its induced
action on $S$; conversely, such a lift preserves original vertices
and midpoints and therefore descends through suppression of the
midpoints. These operations identify the two lift extensions.
We continue to use the induced action of $\operatorname{Aut}(G)$, rather than replacing
$\operatorname{Aut}(G)$ by the possibly larger group $\operatorname{Aut}(S)$.

Fix an original vertex $r\in V(G)$ as a root. We also write $r$ for
its basis vector in $C_0(S;\mathbb Z)$. The following is the integral
chain version of Baez's path description of covering automorphisms
\cite[Lemma~10]{Baez2026}.

\begin{lemma}
\label{lem:baez-chain-lifts}
The group $\operatorname{Cov}(G)$ is isomorphic, as an extension of
$\operatorname{Aut}(G)$ by $H_1(G;\mathbb Z)$, to
\[
 \mathcal L(G)=\{(h,c)\in \operatorname{Aut}(G)\times C_1(S;\mathbb Z):
                        \partial c=hr-r\},
\]
with multiplication
\begin{equation}
 (h,c)(k,d)=(hk,c+hd).
 \label{eq:baez-chain-product}
\end{equation}
Here $hk=h\circ k$, the chain $hd$ is obtained from $d$ by the
linear action just defined, and $c+hd$ is addition in
$C_1(S;\mathbb Z)$.
In particular, a section of $\rho$ is equivalent to a family of
integral chains $(c_h)_{h\in \operatorname{Aut}(G)}$ satisfying
\begin{align}
 \partial c_h&=hr-r,\label{eq:baez-chain-boundary}\\
 c_{hk}&=c_h+h c_k\qquad(h,k\in \operatorname{Aut}(G)).
 \label{eq:baez-chain-cocycle}
\end{align}
\end{lemma}

\begin{proof}
In the path model of the maximal abelian cover, vertices over $v$
are represented by walks from $r$ to $v$. Two such walks represent
the same vertex exactly when their difference has zero integral
homology class. Since a graph has no $2$-chains, this is equivalent
to equality of their signed $1$-chains.

Every integral chain $c$ with $\partial c=hr-r$ is the chain of
some walk from $r$ to $hr$. Indeed, choose one such walk $P$.
Then $c-c(P)$ is an integral cycle. The homomorphism
$\pi_1(S,r)\to H_1(S;\mathbb Z)$ is surjective, so this cycle is
represented by a closed walk $W$ based at $r$. The concatenation
$WP$ has chain $c$.

A lift of $h$ is uniquely determined by the image of the chosen
vertex over $r$, and this image is represented by a walk $Q_h$
from $r$ to $hr$. On the path model, that lift sends
\[
 [P]\longmapsto[Q_h\,h(P)].
\]
Thus its chain coordinate is $c(Q_h)$. Composing the lifts
represented by $Q_h$ and $Q_k$ gives the walk $Q_hh(Q_k)$,
whose chain is $c(Q_h)+h c(Q_k)$. This proves
\eqref{eq:baez-chain-product}. The identity is $(1,0)$ and the
inverse of $(h,c)$ is $(h^{-1},-h^{-1}c)$.

The kernel consists of $(1,z)$ with $\partial z=0$, and hence is
$H_1(S;\mathbb Z)\cong H_1(G;\mathbb Z)$. Conjugation by $(h,c)$ sends $(1,z)$ to
$(1,hz)$, giving the natural homology action. Finally, the map
$h\mapsto(h,c_h)$ is a section precisely when
\eqref{eq:baez-chain-boundary} and \eqref{eq:baez-chain-cocycle}
hold.
\end{proof}

The next lemma concerns a finite group acting by coordinate
permutations. It shows that every family of integral vectors
satisfying the displayed compatibility relation can be written
as the differences $hb-b$ of a single integral vector
$b$.\footnote{In the terminology of group cohomology, the
relation $c_{hk}=c_h+h c_k$ is the \emph{$1$-cocycle identity},
and the lemma states that $H^1(\Gamma,\mathbb Z^\Omega)=0$ for
a finite group permuting a finite basis. This follows from
Shapiro's lemma \cite[Chapter~III, Section~6]{Brown1982};
the same vanishing for divisor permutation modules is used in
\cite[Section~5, Proposition~1]{JoynerKsirMelles2012}.}
The proof below uses only averaging and coordinatewise rounding.

\begin{lemma}
\label{lem:baez-permutation-cocycle}
Let a finite group $\Gamma$ act on $\mathbb Z^\Omega$ by permuting its
standard basis. If $(c_h)_{h\in\Gamma}\subseteq\mathbb Z^\Omega$ satisfies
\[
 c_{hk}=c_h+h c_k\qquad(\forall\;h,k\in \Gamma),
\]
then there is $b\in\mathbb Z^\Omega$ such that $c_h=hb-b$ for every
$h\in \Gamma$.
\end{lemma}

\begin{proof}
First work over $\mathbb Q$ and put
\[
 a:=-\frac1{|\Gamma|}\sum_{k\in \Gamma}c_k.
\]
Rearranging the identity in the statement gives
$h c_k=c_{hk}-c_h$. Reindexing the sum by $k\mapsto hk$ therefore gives
\[
 ha=-\frac1{|\Gamma|}\sum_{k\in \Gamma}(c_{hk}-c_h)=a+c_h.
\]
Let $b:=\lfloor a\rfloor$, taking the floor of each coordinate.
Coordinate permutations commute with this operation, and adding
an integral vector commutes with it as well. Hence
\[
 hb=\lfloor ha\rfloor=\lfloor a+c_h\rfloor
    =\lfloor a\rfloor+c_h=b+c_h.
\]
Thus $c_h=hb-b$, with $b$ integral as required.
\end{proof}

For an integral $0$-chain $Z=\sum_{v\in V(S)}n_vv$, define its
\emph{degree} by $\deg Z=\sum_v n_v$. The coefficients $n_v$ may
be negative. This freedom is essential: an invariant integral
chain of degree one need not be a single fixed vertex.

\begin{proposition}
\label{prop:baez-degree-one}
The extension \eqref{eq:baez-lift-extension} splits if and only if
there exists $Z\in C_0(S;\mathbb Z)$ such that
\[
 hZ=Z\quad(h\in \operatorname{Aut}(G)),\qquad \deg Z=1.
\]
\end{proposition}

\begin{proof}
Suppose first that a section exists. By
Lemma~\ref{lem:baez-chain-lifts}, it supplies chains $c_h$
satisfying \eqref{eq:baez-chain-boundary}--\eqref{eq:baez-chain-cocycle}.
Since $\operatorname{Aut}(G)$ permutes the oriented half-edge basis of $C_1(S;\mathbb Z)$,
Lemma~\ref{lem:baez-permutation-cocycle} gives an integral chain $b$
such that $c_h=hb-b$. Set
\[
 Z=r-\partial b.
\]
Every boundary has degree zero, so $\deg Z=1$. Moreover,
\[
 hZ-Z=hr-r-\partial(hb-b)=hr-r-\partial c_h=0.
\]
Thus $Z$ is invariant.

Conversely, suppose that $Z$ is invariant and has degree one.
Because $S$ is connected, every integral $0$-chain of degree zero
is an integral boundary. To see this explicitly, choose, for each
$v\ne r$, a walk $P_v$ from $r$ to $v$. If
$Y=\sum_v n_vv$ has degree zero, then $n_r+\sum_{v\neq r} n_v=0\implies n_r=-\sum_{v\neq r}n_v$ and
\[
 Y=\sum_{v\ne r}n_v(v-r)
   =\partial\Bigl(\sum_{v\ne r}n_vc(P_v)\Bigr).
\]
Applying this to $Y=r-Z$ gives $b\in C_1(S;\mathbb Z)$ with
$\partial b=r-Z$. Define $c_h=hb-b$. Invariance of $Z$ gives
\[
 \partial c_h=h(r-Z)-(r-Z)=hr-r,
\]
and direct calculation gives
\[
 c_{hk}=hkb-b=(hb-b)+h(kb-b)=c_h+h c_k.
\]
Lemma~\ref{lem:baez-chain-lifts} now supplies a section.
\end{proof}

\begin{proof}[Proof of Theorem~\ref{thm:baez-orbit-splitting}]
\phantomsection
\label{proof:baez-orbit-splitting}

The vertices of $S$ are the original vertices of $G$ and the
midpoints $m_e$ of its unoriented edges. Their $\operatorname{Aut}(G)$-orbit sizes are
therefore precisely the integers appearing in
\eqref{eq:baez-orbit-gcd}.

An integral $0$-chain is $\operatorname{Aut}(G)$-invariant exactly when its coefficients
are constant on each orbit. Consequently, every invariant chain
is an integral linear combination of the orbit sums
\[
 \sum_{v\in V_i}v
 \quad\text{and}\quad
 \sum_{e\in E_j}m_e.
\]
Their degrees are $|V_i|$ and $|E_j|$, respectively. By B\'ezout's
identity, the set of degrees of invariant integral chains is
exactly $d_{\mathrm{orb}}(G)\mathbb Z$. Proposition~\ref{prop:baez-degree-one}
therefore shows that a section exists exactly when $d_{\mathrm{orb}}(G)=1$.

For a section $s$, every element of $\operatorname{Cov}(G)$ is
uniquely of the form $\iota(a)s(h)$. The conjugation calculation
in Lemma~\ref{lem:baez-chain-lifts} shows that multiplication in
these coordinates is
\[
 (a,h)(a',h')=(a+h a',hh'),
\]
which proves the asserted semidirect-product description.
\end{proof}

The next corollary is a formal consequence of
Theorem~\ref{thm:baez-orbit-splitting}. Its abstract-group
formulation uses the torsion-freeness of the deck group and the
standard characterization of split extensions by complements
\cite[Section~4, p.~122]{MalnicPozar2016}.

\begin{corollary}
\label{cor:baez-semidirect}
Let $G$ be a finite connected graph. Then
\[
  \operatorname{Cov}(G)
  \cong H_1(G;\mathbb Z)\rtimes\operatorname{Aut}(G)
  \quad\Longleftrightarrow\quad
  d_{\mathrm{orb}}(G)=1,
\]
where the semidirect product uses the natural action on
integral homology. The equivalence holds even when the
displayed isomorphism is understood only as an abstract
group isomorphism.
\end{corollary}

\begin{proof}
We first record an elementary group-theoretic observation.
For an exact sequence
\[
  1\longrightarrow N\longrightarrow E
    \xrightarrow{q}Q\longrightarrow 1
\]
with $N$ torsion-free and $Q$ finite, the extension splits if
and only if $E$ contains a subgroup $K$ of order $|Q|$.
Indeed, a section supplies such a subgroup. Conversely,
$K\cap N=\{1\}$, since $K$ is finite and $N$ is torsion-free.
Thus $q|_K$ is injective and, by equality of orders, is an
isomorphism onto $Q$. Its inverse, followed by inclusion into
$E$, is a section.

Here the kernel $H_1(G;\mathbb Z)$ is free abelian and
$\operatorname{Aut}(G)$ is finite. If the displayed abstract
isomorphism exists, it takes the finite factor
$\{0\}\times\operatorname{Aut}(G)$ to a subgroup of
$\operatorname{Cov}(G)$ of order $|\operatorname{Aut}(G)|$.
The observation above therefore shows that the canonical
extension splits, and Theorem~\ref{thm:baez-orbit-splitting}
gives $d_{\mathrm{orb}}(G)=1$.

Conversely, if $d_{\mathrm{orb}}(G)=1$, the same theorem gives
the stated semidirect-product isomorphism, with the natural
homology action.
\end{proof}

If $\operatorname{Aut}(G)$ fixes a vertex of $G$, then one vertex orbit has size one,
so the extension splits. This recovers Baez's fixed-basepoint
sufficient condition \cite[Section~8]{Baez2026}. The same conclusion
holds if $\operatorname{Aut}(G)$ stabilizes an unoriented edge, even if some elements
reverse it: its midpoint is fixed in $S(G)$. Neither condition is
necessary, as Example~\ref{ex:baez-k23} below shows.

\begin{corollary}
\label{cor:baez-sylow}
The extension \eqref{eq:baez-lift-extension} splits if and only if,
for every prime $\ell$ dividing $|\operatorname{Aut}(G)|$, a Sylow $\ell$-subgroup of
$\operatorname{Aut}(G)$ fixes a vertex of $G$ or stabilizes an unoriented edge setwise.
For each prime, ``a Sylow subgroup'' may equivalently be replaced
by ``every Sylow subgroup''; the fixed vertex or edge may vary.
\end{corollary}

\begin{proof}
If $x$ is a vertex or an unoriented edge and $\operatorname{Aut}(G)_x$ is its stabilizer,
then $|\operatorname{Aut}(G)\cdot x|=[\operatorname{Aut}(G):\operatorname{Aut}(G)_x]$. This index is prime to $\ell$ exactly
when $\operatorname{Aut}(G)_x$ contains a Sylow $\ell$-subgroup of $\operatorname{Aut}(G)$. All Sylow
$\ell$-subgroups are conjugate, so a fixed object for one gives a
fixed object for each conjugate. Every orbit size divides $|\operatorname{Aut}(G)|$.
Thus their gcd is one exactly when, for each prime divisor $\ell$
of $|\operatorname{Aut}(G)|$, at least one orbit size is prime to $\ell$. Apply
Theorem~\ref{thm:baez-orbit-splitting}.
\end{proof}

The same proof applies to any subgroup $K\leq\operatorname{Aut}(G)$.
The restricted extension
$$1\longrightarrow H_1(G;\mathbb Z)\longrightarrow\rho^{-1}(K)\longrightarrow K\longrightarrow 1$$
splits exactly when the gcd of
the sizes of the $K$-orbits on vertices and unoriented edges is one.

\begin{example}[Dipoles]
\label{ex:baez-dipoles}
Let $D_m$ be the graph with two vertices $u,v$ joined by $m\geq1$
parallel edges. Its automorphism group is $S_m\times C_2$:
$S_m$ permutes the edges and $C_2$ interchanges the two vertices.
There is one vertex orbit of size two and one unoriented-edge
orbit of size $m$. Hence
\[
 d_{\mathrm{orb}}(D_m)=\gcd(2,m),
\]
and the canonical lift extension splits if and only if $m$ is odd.

When $m=2k+1$, the invariant degree-one chain can be written
explicitly. If $m_1,\ldots,m_m$ are the edge midpoints, take
\[
 Z=\sum_{i=1}^{m}m_i-k(u+v).
\]
Its degree is $m-2k=1$, so
Proposition~\ref{prop:baez-degree-one} constructs a section.
For even $m$, every invariant chain on $S(D_m)$ has even degree,
which obstructs a section.

For $m=3$, the maximal abelian cover is the hexagonal net;
for $m=4$, it is the diamond net \cite[Section 8]{Baez2026}.
Their canonical extensions split and do not split, respectively.
Baez already identifies the symmetries fixing a hexagon centre
\cite[Section~8]{Baez2026}. This order-$12$ subgroup
is a complement to the canonical deck group: it intersects that
torsion-free group trivially and maps isomorphically onto
$\operatorname{Aut}(D_3)$, which also has order $12$.
Thus the nonsplitting assertion in the same section needs
correction. The invariant chain $m_1+m_2+m_3-u-v$ gives an
algebraic verification.
\end{example}

\begin{example}[A split extension without a global fixed point]
\label{ex:baez-k23}
Let $G=K_{2,3}$, with bipartition $U\sqcup W$, where $|U|=2$ and
$|W|=3$. Its automorphism group is $S_2\times S_3$. The two vertex
orbits have sizes two and three, and the unique unoriented-edge
orbit has size six. Therefore
\[
 d_{\mathrm{orb}}(K_{2,3})=\gcd(2,3,6)=1.
\]
The invariant integral chain
\[
 Z=\sum_{w\in W}w-\sum_{u\in U}u
\]
has degree one. Hence the extension splits, although the full
automorphism group fixes no vertex and stabilizes no individual
edge. This also illustrates why integral chains with negative
coefficients are needed in Proposition~\ref{prop:baez-degree-one}.
\end{example}

The following corollary collects the complete answers to both
questions posed by Baez, already established in
Corollary~\ref{cor:abelian-cover-descent} and
Theorem~\ref{thm:baez-orbit-splitting}.
A rooted tree is \emph{root-rigid} if its only automorphism
fixing the root is the identity.

\begin{corollary}
\label{cor:baez-combined}
Let $G$ be a finite connected graph, and let $C:=\Core(G)$.
When $C\neq\varnothing$, write $G$ as the graph obtained from
$C$ by attaching a finite rooted tree $T_v$ at each
$v\in V(C)$, allowing $T_v$ to consist only of its root. Then:
\begin{enumerate}
 \item Every automorphism of $G^{\mathrm{ab}}$ descends to $G$
       if and only if either $G$ is a tree, or
       $C\neq\varnothing$ and every $T_v$ is root-rigid.
 \item The canonical extension
       \eqref{eq:baez-lift-extension} splits if and only if
       $d_{\mathrm{orb}}(G)=1$.
\end{enumerate}
If $G$ is a tree, then $G^{\mathrm{ab}}\cong G$ and
$H_1(G;\mathbb Z)=0$, so both descent and splitting are automatic.
Whenever both properties hold,
\[
 \operatorname{Aut}(G^{\mathrm{ab}})
   \cong H_1(G;\mathbb Z)\rtimes\operatorname{Aut}(G),
\]
with the natural homology action, the canonical deck subgroup,
and the canonical descent map.
\end{corollary}

\section{Conclusion}\label{sec:scope}

The maximal abelian cover determines a finite connected base graph
unless the base is unicyclic. The proof separates two pieces of
information: bounded-displacement automorphisms recover the core
covering, and the finite trees attached to that core recover the
remaining vertices. For a unicyclic base, only the primitive pattern
of rooted-tree decorations survives.

Section~\ref{sec:baez} gives complete criteria for both questions
posed by Baez. For a graph with nonempty core, every automorphism
of its maximal abelian cover descends if and only if each rooted
tree attached to the core is root-rigid; for a tree, descent is
automatic (Corollary~\ref{cor:abelian-cover-descent}).

The canonical lift extension splits if and only if the greatest
common divisor of the sizes of the $\operatorname{Aut}(G)$-orbits
on vertices and unoriented edges is one
(Theorem~\ref{thm:baez-orbit-splitting}). Equivalently, every Sylow
subgroup fixes a vertex or stabilizes an unoriented edge
(Corollary~\ref{cor:baez-sylow}). When both descent and splitting
hold, these criteria give
\[
  \operatorname{Aut}(G^{\mathrm{ab}})
  \cong H_1(G;\mathbb Z)\rtimes\operatorname{Aut}(G),
\]
with the natural homology action, the canonical deck subgroup,
and the canonical descent map
(Corollary~\ref{cor:baez-combined}).

One may ask which other characteristic quotients of the
fundamental group retain enough information to reconstruct the base.
For abelian covers that are not maximal, the rank of the deck group
can be smaller than $\beta(G)$, and the final rank comparison in
Theorem~\ref{thm:deck} no longer forces the intermediate covering
to have one sheet. Identifying useful rigidity hypotheses in this
setting remains a natural direction.

Finally, having isomorphic maximal abelian covers is stronger than
having \emph{some} common connected regular cover with abelian deck
groups. For instance, let $G$ have cycle rank at least two and let
$Q\to G$ be a connected regular cyclic covering of degree $d>1$.
Such a cover is obtained by assigning voltage $1\in\Z/d\Z$ to one
edge outside a spanning tree and $0$ to the other chosen edges.
Then $Q$ is a common abelian cover of $G$ and $Q$, through this
projection and the identity map. However,
\[
  \beta(Q)-1=d\bigl(\beta(G)-1\bigr),
\]
so their maximal abelian covers are not isomorphic by
Theorem~\ref{thm:classification}. A characterization of common
abelian covers must therefore address a different comparison problem.

\subsubsection*{Funding information}

The authors acknowledge support from Conselho Nacional de Desenvolvimento Científico e Tecnológico (CNPq), and Fundação de Amparo à Pesquisa do Estado de Minas Gerais (FAPEMIG).

\subsubsection*{Use of AI}

Authors acknowledge the use of generative AI (GPT-6) for the preparation of this manuscript. Authors remain solely responsible for the content and accuracy of the work.

\bibliographystyle{plain}
\bibliography{References}

\end{document}